\documentclass[12pt,reqno]{amsart}

\usepackage{epsfig}
\usepackage{amscd}
\usepackage[mathscr]{eucal}
\usepackage{amssymb}
\usepackage{amsxtra}
\usepackage{amsmath}
\usepackage[all]{xy}
\usepackage{mathtools}
\usepackage[pdfencoding=auto]{hyperref}
\usepackage{bookmark}

\DeclarePairedDelimiter\abs{\lvert}{\rvert}

\theoremstyle{plain}

\newtheorem{prop}{Proposition}[section]
\newtheorem{cor}[prop]{Corollary}

\usepackage{amsthm}
\newtheorem{theorem}{Theorem}[section]

\theoremstyle{definition}

\theoremstyle{remark}

\begin{document}

\bigskip

\title[An approx. fun. equ. with exponentially decaying error]{An approximate functional equation for the Riemann zeta function with exponentially decaying error}

\date{\today}

\author{Yochay Jerby}

\address{Yochay Jerby, Faculty of Sciences, Holon Institute of Technology, Holon, 5810201, Israel}

%
%
% ...

\begin{abstract} It is known by a formula of Hasse-Sondow that the Riemann zeta function is given, for any $ s=\sigma+it \in \mathbb{C}$, by $ \sum_{n=0}^{\infty} \widetilde{A}(n,s)$ where $$ 
\widetilde{A}(n,s):=\frac{1}{2^{n+1}(1-2^{1-s})} \sum_{k=0}^n \binom{n}{k}\frac{(-1)^k}{(k+1)^s} .$$ 
We prove the following approximate functional equation for the Hasse-Sondow presentation: For $  \vert t \vert = \pi xy $ and $  2y \neq (2N-1)\pi $ then 
$$ 
\zeta(s)= \sum_{n \leq x } \widetilde{A}(n,s)+\frac{\chi(s)}{1-2^{s-1}} \left (\sum_{k \leq y} (2k-1)^{s-1} \right ) +O \left (e^{-\omega(x,y) t} \right ), $$
where $ 0 <\omega(x,y)$ is a certain transcendental number determined by $ x$ and $  y$. A central feature of our new approximate functional equation is that its error term is of exponential rate of decay. The proof is based on a study, via saddle point techniques, of the asymptotic properties of the function $$  
\widetilde{A}(u,s):= \frac{1}{2^{u+1} (1-2^{1-s}) \Gamma(s)}    \int_{0}^{\infty}   \left ( e^{-w} \left ( 1- e^{-w} \right)^u \right ) w^{s-1} dw,$$ 
and integrals related to it. \end{abstract}

\maketitle
%
%
\section{Introduction}

The classical approximate functional  equation for the Riemann zeta function was proved by Hardy and Littlewood in the series of works \cite{HL,HL2,HL3}. Let  $A(n,s):=n^{-s}$ and $s=\sigma +it$ with $0  \leq \sigma \leq 1$. The theorem states that the following holds 
\begin{equation}
\label{eq:HL-AFE} 
\zeta(s)= \sum_{n \leq x } A(n,s) + \chi(s)  \left ( \sum_{k \leq y } k^{s-1} \right )  + O \left ( x^{-\sigma} + x^{\frac{1}{2} - \sigma} y^{-\frac{1}{2}}  \right ),
\end{equation} when $x,y \geq 1$ are such that $\vert t \vert = 2 \pi x y$ and 
$\chi(s):=2^s \pi^{s-1}sin \left ( \frac{\pi s}{2}\right ) \Gamma(1-s)$  is the function appearing in the functional equation $\zeta(s) = \chi(s) \zeta(1-s)$. An asymptotic expansion for the error term in (\ref{eq:HL-AFE}), is further described by the celebrated Riemann-Siegel formula, published in 1932  by C. L. Siegel, based on original manuscripts of Riemann, see \cite{Si} and also \cite{E,I,T}.  

In 1994, Sondow proved, via Euler's transformation of series, that the Riemann zeta function can be expressed as  
 \begin{equation} 
  \label{eq:HS1}
  \zeta(s) := \sum_{n=0}^{\infty} \widetilde{A}(n,s) 
  \end{equation}
   where 
 \begin{equation} 
 \label{eq:HS2}
 \widetilde{A}(n,s):=\frac{1}{2^{n+1}(1-2^{1-s})} \sum_{k=0}^n \binom{n}{k}\frac{(-1)^k}{(k+1)^{s}},
\end{equation}  
for any $s \in \mathbb{C}$, see \cite{SE}. Formula (\ref{eq:HS1}) also appears in the appendix to a work by Hasse from 1930, where he proved results of similar nature, see \cite{H}. In the appendix, Hasse attributes the conjecture of formula (\ref{eq:HS1}) to Knopp. In this work we prove the following analog of the classical approximate functional equation, for the Hasse-Sondow representation of zeta:
\begin{theorem}
\label{thm:main-theorem} 
 For $x,y \geq 1$ such that $\vert t \vert = \pi xy $ and $2y \neq (2N-1)\pi $ the following holds 
\begin{equation}
\label{eq:NEW-AFE}
\zeta(s)= \sum_{n \leq x } \widetilde{A}(n,s) +\frac{\chi(s)}{1-2^{s-1}} \left (\sum_{k \leq y} (2k-1)^{s-1} \right ) +O \left (e^{-\omega(x,y) t} \right ), 
\end{equation} where $\omega(x,y)>0$ is a given transcendental number, uniquely determined by $x$ and $y$. 
\end{theorem}
  
A central feature of our new approximate functional equation (\ref{eq:NEW-AFE}) is that its error term, $O \left (e^{-\omega(x,y) t} \right )$, is of exponential rate of decay, in contrast to the fractional power rate of decay in the classical case. Moreover, already the main term of our new equation (\ref{eq:NEW-AFE}) gives an approximation of $\zeta(s)$ which attains the level of accuracy of the asymptotic Riemann-Siegel formula, expanded to its conjectured optimal order\footnote{In \cite{Be} Berry conjectured, supported by extensive numerical verification, that the best possible accuracy of the Riemann-Siegel formula is $O \left ( e^{-\pi t} \right )$, obtained by expanding the formula to order $r^{\ast}(t)= \left [ 2 \pi t \right]$. That is, by adding the first $r^{\ast}(t)$ correcting terms to the main sum. For instance, for the first zero of zeta, the optimal accuracy requires adding the first $r^{\ast}(t)=89$ correcting terms to the main sum. In particular, as $t$ grows the number of correcting terms required to obtain exponential accuracy becomes infinite, see also \cite{BK}.}. The exact definition of $\omega(x,y)$ is given below.   

The proof of Theorem \ref{thm:main-theorem}, as the proof of the classical approximate functional equation, (\ref{eq:HL-AFE}), in \cite{HL}, utilizes  the Euler-Maclaurin summation formula of order one which states that 
\begin{equation}
\label{eq:EM}
\sum_{n=M}^N A(n,s) = \int_M^N A(u,s) du + \int_M^N \psi(u) \left ( \frac{\partial}{\partial u} A(u,s) \right ) du+\frac{1}{2} \left [ A(M,s)+A(N,s) \right ],
\end{equation} 
where $A(u,s)$ is a continuous function in the $u$-variable and
$
\psi(u):=u - [u]-\frac{1}{2}$.
In the classical Hardy-Littlewood case\footnote{ Hardy and Littlewood gave two proofs of the approximate functional equation (\ref{eq:HL-AFE}). The "first proof" is based on real analysis and the Eular-Maclaurin summation formula and is presented in \cite{HL,HL3}  while the "second proof" is based on complex analysis and contour integration and appears in \cite{HL2}. The Riemann-Siegel formula is obtained by further applying saddle point techniques to the integrals of \cite{HL2}. In our case, of formula (\ref{eq:NEW-AFE}), the interpolating function $\widetilde{A}(u,s)$, of the Euler-Mclaurain summation 
itself, is represented by an integral whose asymptotic analysis requires application 
of the saddle-point technique. Therefore, although our starting point is the Euler-Mclaurain formula (as in the "first proof" of \cite{HL}), the proof of Theorem \ref{thm:main-theorem} eventually involves a combination the various methods.} of (\ref{eq:HL-AFE}), the function $A(u,s)=u^{-s}$ is readily taken. However, in our case of the Hasse-Sondow presentation, the required interpolating function turns to be
\begin{equation}
\label{eq:GL-integral}
\widetilde{A}(u,s):= \frac{1}{2^{u+1} (1-2^{1-s}) \Gamma(s)}  \int_{0}^{\infty}   \left ( e^{-w} \left ( 1- e^{-w} \right)^u \right ) w^{s-1} dw.
\end{equation}

Let us note that $\widetilde{A}(u,s)$ is related to, $\Phi^{\ast}_u(z,s,a)$, the generalized Hurwitz-Lerch zeta function of Goyal and Laddha, via $
\widetilde{A}(u,s)=\frac{\Phi^{\ast}_{-u}(1,s,1)}{2^{u+1}(1-2^{1-s})}$, see \cite{GL} and also \cite{GN,NC} in relations to nuclear and condensed matter physics. In Section 2, we include a self-contained derivation of (\ref{eq:GL-integral}) from the corresponding N{\o}rlund-Rice integral, see \cite{FS}.

 Applying the Euler-Maclaurin formula (\ref{eq:EM}) to $\widetilde{A}(u,s)$, gives
\begin{equation}
\label{eq:EM2}
\zeta(s)-\sum_{n \leq x} \widetilde{A}(u,s) = \int_x^{\infty} \widetilde{A}(u,s) du + \int_x^{\infty} \psi(u) \left ( \frac{\partial}{\partial u} \widetilde{A}(u,s) \right ) du+\frac{1}{2}  \widetilde{A}(x,s).
\end{equation} 

The proof of Theorem \ref{thm:main-theorem}, thus, relies on the asymptotic saddle-point approximation of the three integrals appearing in the right hand side of (\ref{eq:EM2}). That is:

\bigskip

\begin{enumerate}

\item The function $\widetilde{A}(x,s)$ itself. In Corollary \ref{cor:2.3} we show that $\widetilde{A}(u,s)=O \left ( e^{- \omega(x,y) t} \right).$ 

\bigskip

\item  In Theorem \ref{thm:first-integral} we show that 
$
\widetilde{I}_1(x,s):= \frac{\chi(s)}{2^{s-1}-1} \left (\sum_{k \leq y} (2k-1)^{s-1} \right ) +O \left (e^{-\omega(x,y) t} \right )$. In particular, the main sum arises from pole contributions. 

\bigskip

\item In Theorem \ref{thm:int2} we show that also $\widetilde{I}_2(x,s)=O \left ( e^{- \omega(x,y) t} \right)$.
\bigskip

\end{enumerate}

A key feature is that the saddle points of all three integrals (1)-(3) are shown to be given as the solutions of the following transcendental equation $
e^{w}-\alpha \cdot i w -1=0$. The solutions $
w_k(\alpha)=x_k(\alpha)+y_k(\alpha) \cdot i$ can be ordered by $k \in \mathbb{Z}$ according to the value of $y_k(\alpha)$. Moreover, the main contribution for the asymptotic behavior of the integrals is shown to come from the $k(\alpha)$-th saddle point, which is the saddle point whose $\vert x_{k(\alpha)}(\alpha) \vert$ value is minimal. Finally, the exponential decay rate of Theorem \ref{thm:main-theorem} is given, in Corollary \ref{cor:2.3}, as follows \begin{equation} 
\omega(x,y):=- \alpha \cdot log \abs*{\frac{1-e^{-w_{k(\alpha)}(\alpha)}}{2}} + \theta_{k(\alpha)}(\alpha)-\frac{\pi}{2}, 
\end{equation}
where 
$
\alpha=\alpha(x,y):=\frac{t}{x}=\frac{1}{\pi y}$,
and $w_k(\alpha):=r_k(\alpha)e^{i \theta_k(\alpha)}$ is the polar presentation of the $k(\alpha)$-th saddle-point. In particular, the reason for excluding the case $2y \neq (2N-1)\pi $ in Theorem \ref{thm:main-theorem} is that these are the degenerate values for which $\omega(x,y)=0$.

\bigskip

The rest of the work is organized as follows: In Section 2 we review the definition of $\widetilde{A}(u,s)$ and study its asymptotic properties. The integral $\widetilde{I}_1(x,s)$ is investigated in Section 3. The integral $\widetilde{I}_2(x,s)$ is considered in section 4 together with the conclusion of the proof of the main theorem. In Section 5 we present concluding remarks and discuss possible relations to quantum chaos. 
\section{The Function \texorpdfstring{$\widetilde{A}(u,s)$}{$A(u,s)$} and its Asymptotic Properties}
\label{sec:2}

 Recall that the $n$-th difference of a continuous function $f(u)$ is given by the sum 
\begin{equation}
 D_n [f]:= \sum_{k=0}^n \binom{n}{k} (-1)^k f(k). 
 \end{equation} For instance, the $n$-th coefficient of the Hasse-Sondow presentation (\ref{eq:HS2}), could be written as 
\begin{equation} 
\label{eq:def}
\widetilde{A}(n,s)=\frac{D_n[A(v+1,s)]}{2^{n+1}(1-2^{1-s})},
\end{equation} where $A(v,s)=v^{-s}$. 

In general, one observes that $n$-th differences satisfy the following trivial "naive" bound $
\vert D_n [f] \vert \leq 2^{n} \max_{0 \leq k \leq n} \vert f(k) \vert$, which implies, assuming that $f(u)$ is bounded, that $D_n[f]=O(2^n)$. In the mid 1960's, the study of asymptotic properties of $n$-th differences, $D_n[f]$, was found to be of central importance for various questions in algorithm analysis related to the works of De Bruijn, Knuth and Rice, see for instance Section 5.2.2 and 6.3 of \cite{K}. As a result, various saddle-point techniques and countour integration methods were developed and applied for the study of such asymptotics. In particular, one of the central phenomena, observed by the algorithm analysis community, is that $n$-differences actually tend to exhibit asymptotic behavior that is essentially slower than $O(2^n)$, a phenomena known as exponential cancellation, see \cite{FS}. 

As mentioned in the introduction, in order to apply analytic methods for the investigation of the asymptotic behavior of $D_n[f]$, a continuous function, interpolating the values of $D_n[f]$ is required. In general, such a function is given via the N{\o}rlund-Rice integral as follows: If $f(v)$ is an analytic function in a domain containing $0<Re(s)$, one has 
\begin{equation} Res_{v=k} f(v) \frac{n!}{v(v-1)...(v-n)}=(-1)^{n-k} \binom{n}{k} f(k),
\end{equation} for $0 \leq k \leq n$. Note that 
\begin{equation} 
\frac{n!}{v(v-1)...(v-n)}= (-1)^{k-1} B(n+1,-v),
\end{equation} where $B(z,w)$ is the beta function. Hence, the $n$-th difference could be expressed in terms of the following N{\o}rlund-Rice integral formula 
\begin{equation} 
D_n[f]=\frac{1}{2 \pi i} \int_C B(n+1,-v)f(v) dv
\end{equation} 
where $C$ is a contour of integration encircling $[0,n]$, avoiding integral values between zero and $n$.
Moreover, under mild conditions, the N{\o}rlund-Rice integral further simplifies, and takes the following form 
\begin{equation} 
\label{eq:NR} D_n[f] :=\frac{1}{2 \pi i} \int_{c-i \infty}^{c+i \infty} B(n+1,-v)f(v) dv, 
\end{equation} 
for $c<0$. In our case we have: 

\begin{prop}
\label{prop:integral}
For $s=\sigma+it$ with $\sigma>1-u$ the required interpolating function of (\ref{eq:def}) is given by 
\begin{equation} 
\label{eq:GL}
\widetilde{A}(u,s):= \frac{1}{2^{u+1} (1-2^{1-s}) \Gamma(s)}  \int_{0}^{\infty}   \left ( e^{-w} \left ( 1- e^{-w} \right)^u \right ) w^{s-1} dw.
\end{equation} 
\end{prop} 

\hspace{-0.6cm} \bf Proof: \rm In view of the above, we can express 
\begin{equation} 
D_n[A(v+1,s)]:=  \frac{1}{2 \pi i} \int_{c-i \infty}^{c+i \infty} B(n+1,-v)(v+1)^{-s} dv.
 \end{equation}
From the definition of the beta function 
$
B(z,w)= \int_0^1  t^{z-1} (1-t)^{w-1}dt$, it follows that, for $Re(z),Re(w)>0$, the following holds 
\begin{multline}
\label{eq:Mellin1}
 D_n[A(v+1,s)] =\frac{1}{2 \pi i} \int_{c-i \infty}^{c+i \infty}  B(n+1,-v)(v+1)^{-s} dv = \\ = \frac{1}{2 \pi i} \int_{c-i \infty}^{c+i \infty} \left ( \int_0^1  t^n (1-t)^{-(v+1)} dt \right ) (v+1)^{-s} dv   
=  \\=  \int_0^1 t^n \left ( \frac{1}{2 \pi i} \int_{c-i \infty}^{c+i \infty} (1-t)^{-(v+1)} (v+1)^{-s} dv \right ) dt.
\end{multline} Consider the function
 \begin{equation} G(t,s):= \frac{1}{2 \pi i} \int_{c-i \infty}^{c+i \infty} t^{-v} (v+1)^{-s} dv,
 \end{equation} which is, by definition, the inverse Mellin transform of $(v+1)^{-s}$. On the other hand, note that 
 \begin{multline} \int_0^1 \left (ln(t)^{s-1} \right ) t^v dt =\left [\begin{array}{cc} z=ln(t) \\ dz =\frac{dt}{t} \end{array} \right ] = - e^{i \pi s} \int_{0}^{\infty}  z^{s-1}  e^{-(z+1)u} du= \\ =\left [\begin{array}{cc} w=(v+1)z \\ dw=(v+1)dz \end{array} \right ]=-e^{i \pi s} (v+1)^{-s} \int_0^{\infty} w^{s-1} e^{-w} dw =  -e^{i \pi s}\Gamma(s)  (v+1)^{-s}.
 \end{multline}
  Hence, for $Re(s)>0$, we have 
  \begin{equation} G(t,s) := \left \{ 
\begin{array}{cc} -\frac{ln(t)^{s - 1} \cdot t }{e^{ \pi i s}\Gamma(s)} & 0 \leq t \leq 1 \\ 0 & 1<t \end{array} \right. 
\end{equation}   
In particular, taking $c>-1$, equation (\ref{eq:Mellin1}) could be written as 
\begin{multline} 
D_n[A(v+1,s)]
=    \int_0^1   \left (\frac{G(1-t,s)}{1-t} \right ) t^n dt= -\frac{1}{e^{i \pi s} \Gamma(s)}    \int_0^1    \left ( ln(1-t)^{s-1} \right ) t^n dt = \\ =-\frac{1}{e^{i \pi s}\Gamma(s)}    \int_0^1    (1-t)^n \left ( ln(t)^{s-1} \right )  dt=\left [\begin{array}{cc} w=ln(t) \\ dw =\frac{dt}{t} \end{array} \right ] = \\ =-\frac{1}{e^{i \pi s} \Gamma(s)}    \int_{-\infty}^0    \left (e^w(1-e^w)^n \right ) w^{s-1} dw =
\frac{1}{\Gamma(s)}    \int_{0}^{\infty}    \left (e^{-w}\left (1-e^{-w} \right )^n \right ) w^{s-1} dw,
\end{multline} 
as required. $\square$

\bigskip 

Proposition \ref{prop:integral} enables us to apply saddle-point techniques\footnote{The asymptotics of $\sum_{k=1}^n\binom{n}{k} \frac{(-1)^k}{k^s}$ (a variant of our $D_n[(A(v+1,s)]$) was studied by Flajolet and Sedgewick via the Rice method for N{\"o}rlund-Rice integrals, in \cite{FS}. It was shown in \cite{FS} that for fixed $s =\sigma+it \in \mathbb{C}$ non-integer, and $n$ arbitrarly large, one has $$\sum_{k=1}^n\binom{n}{k} \frac{(-1)^k}{k^s}=  \frac{1}{2 \pi i} \int_{\frac{1}{2}-i \infty}^{\frac{1}{2}+i \infty} B(u+1,-v)v^{-s} dv = O\left ((log(u)^{\sigma} \right ).$$ In our case, we are concerned with the asymptotics of $D_n[A(v+1,s)]$ when $n$ is proportional to $t=Im(s)$, for which classical saddle-point methods applied to (\ref{eq:GL}) are  adequate.} in order to describe the asymptotic behvior of $\widetilde{A}(u,s)$. We refer the reader to \cite{BH,DB,Wo} for standard references on the subject. First, we are interested in the asymptotic behavior of the $1$-parametric family of integrals
\begin{equation}
\label{eq:saddle}
I(\sigma,t ; \alpha):=  \int_{0}^{\infty}  
e^{-w} \left ( 1- e^{-w} \right)^{\alpha t}   w^{s-1} dw,
\end{equation}
as $t \rightarrow \infty$. Note that for $s=\sigma+it$ and $u=\alpha \cdot t$ equation (\ref{eq:GL}) can be written as 
\begin{equation}
\label{eq:A}
\widetilde{A}(u,s):=\frac{I(\sigma,t; \alpha)}{2^{u+1}(1-2^{1-s}) \Gamma(s)} 
\end{equation}
Let us express 
\begin{equation}
I(\sigma,t ;\alpha) := \int_0^{\infty} g(\sigma,w) e^{t f(w ; \alpha)} dw
\end{equation} where the two functions are given by  
\begin{equation}
\label{eq:29} 
\begin{array}{ccc} g(\sigma,w):=e^{-w} w^{\sigma -1} 
& ; & f(w ;\alpha):=\alpha \cdot ln \left (1-e^{-w} \right ) + i \cdot ln(w). 
\end{array}
\end{equation}
When expressed in such a form, the saddle point technique implies that the asymptotic properties of the integrals $I(\sigma,t ;\alpha)$ are determined by the location of saddle points of $f(w ;\alpha)$ in the complex plane, that is, the by the location of the solutions of the transcendental equation 
\begin{equation}
\label{eq:equ30} \frac{\partial f}{\partial w} (w ; \alpha) = \frac{ \alpha \cdot e^{-w}}{1-e^{-w}}+\frac{i}{w}=0. 
\end{equation}
Note that the solutions of (\ref{eq:equ30}) are equivalent to the  solutions of 
\begin{equation}
\label{eq:trans2}
 e^w-  \alpha \cdot i w -1=0,
\end{equation} excluding the trivial solution at $w=0$. Set $w=x+yi$ and express  (\ref{eq:trans2}) as the system of two equations
\begin{equation}
\label{eq:sys2}
\begin{array}{ccc}
H^{\alpha}_1(x,y):=e^x cos(y)+ \alpha y-1=0  & ; & H^{\alpha}_2(x,y):=e^x sin(y)- \alpha x=0. \end{array}  
\end{equation}  
Figure 1 illustrates the solutions of (\ref{eq:sys2}) in the region $-2 \leq x \leq 2 $ and $-20 \leq y \leq 20$ for $\alpha=\frac{1}{\pi}$. In particular, the saddle points appear at the intersection of the two curves: 

\begin{figure}[ht!]
	\centering
	\includegraphics[scale=0.17]{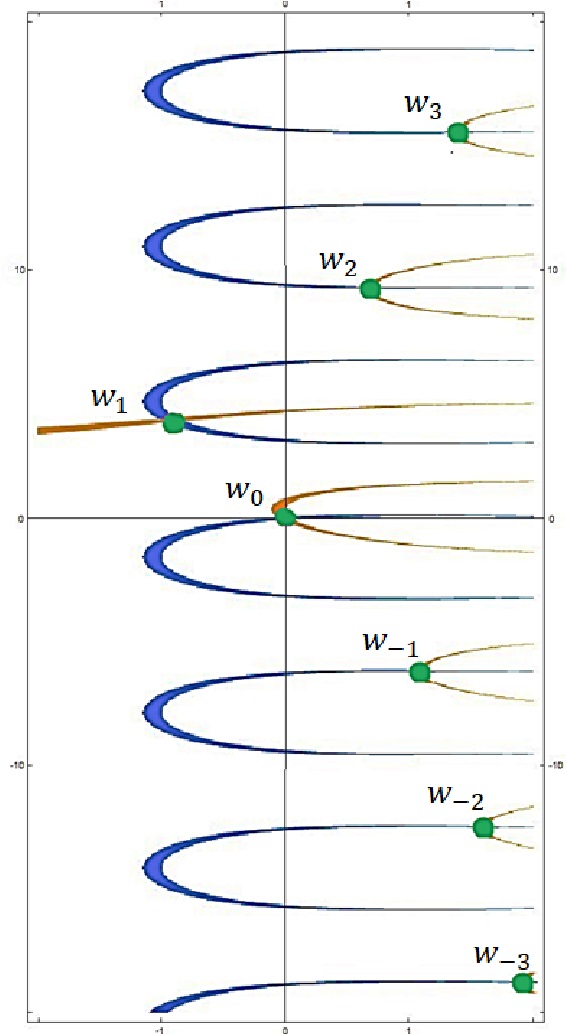}
	\caption{Solutions of $H^{\alpha}_1(x,y)=0$ (brown) and $H^{\alpha}_2(x,y)=0$ (blue) for $-2 \leq x \leq 2,-20 \leq y \leq 20$ and $\alpha=\frac{1}{\pi}$. }
\end{figure} 

We order the solutions $w_k(\alpha)=x_k(\alpha)+y_k(\alpha) i$, with $k \in \mathbb{Z}$, according to the value of $y_k(\alpha)$. As a convention we set $w_0(\alpha)=0$ even though it is a solution only of (\ref{eq:trans2}) and not of (\ref{eq:equ30}). As a direct application of the saddle point technique, we have: 

\begin{prop} 
\label{prop:2.2} The following holds
\begin{equation} \vert 
I(\sigma,t ;\alpha) \vert \sim \abs{ \sqrt{\frac{2 \pi}{t \vert f''(w_{k(\alpha)}(\alpha) ;\alpha ) \vert }}  e^{-w_{k(\alpha)}(\alpha)} \left ( 1-e^{-w_{k(\alpha)}(\alpha)} \right )^{\alpha t} w_{k(\alpha)}(\alpha)^{\sigma +it -1} },
\end{equation}
where for $N>1$ and $\alpha \in  \left (\frac{2}{(2N+1)\pi},\frac{2}{(2N-1) \pi} \right )$
\begin{equation} k(\alpha): = \left \{ \begin{array}{cc} N & \abs{x_N}<\abs{x_{N+1}} \\ N+1 & otherwise \end{array} \right. 
\end{equation} 
and $k(\alpha)=1$ for $\alpha \in \left ( \frac{2}{ \pi} , \infty \right )$. 
\end{prop}

\bigskip

Proposition \ref{prop:2.2} expresses the fact that the main contribution to the asymptotic behavior of $\abs{I(\sigma,t ;\alpha}$ is determined by the saddle point $w_k(\alpha)$ with $\vert x_k(\alpha) \vert$ minimal and $y_k(\alpha)>0$, which is the $k(\alpha)$-th saddle point. For instance, $k(\alpha)=2$ for $\alpha=\frac{1}{\pi}$, as seen in Fig. 1. 

In order to further describe of the general position of the required saddle points, let us consider the graphs of $H_i^{\alpha}(x,y)$, for $i=1,2$, as the value of $x$ varies. Starting with $x<<0$ the two graphs are essentially linear. As the value of $x$ grows, the exponent and trigonometric factors begin to contribute, as illustrated in the following Fig. 2:  

\begin{figure}[ht!]
	\centering
	\includegraphics[scale=0.5]{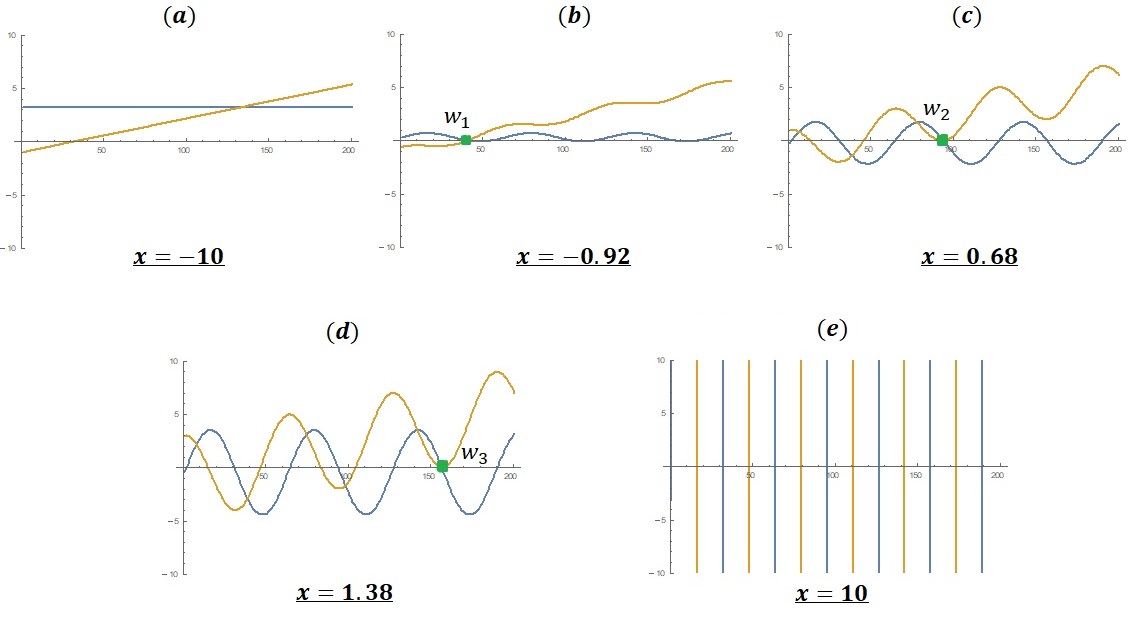}
	\caption{Graphs of $H^{\alpha}_1(x_0,y)$ (brown) and $H^{\alpha}_2(x_0,y)$ (blue) for $0 \leq y \leq 20$ with $\alpha=\frac{1}{\pi}$ and $ x_0=-10 (a),-0.92 (b),0.68 (c),1.38 (d),10 (e)$. }
\end{figure} 

Assume $\alpha \in \left (\frac{2}{(2N+1)\pi},\frac{2}{(2N-1) \pi} \right )$, as in Proposition \ref{prop:2.2}. From the above observation, we see that saddle points below the $\left [\frac{N}{2} \right ]$-th point are to occur around intersections of local maxima of $H^{\alpha}_1(x,y)$ with the axis which, in turn, correspond to local maxima of $cos(y)$. That is, we have $
y_k(\alpha) \approx  2k \pi$ for $k<\frac{N}{2}$. 

For the higher saddle points, the parity of $N$ needs to be considered. For $N=2M$ even, 
we have 
\begin{equation}
\begin{array}{ccc}
y_M(\alpha) \approx \left (2M-\frac{1}{4} \right ) \pi & ; & y_{M+1}(\alpha) \approx \left (2M+\frac{5}{4} \right ) \pi. 
\end{array}
\end{equation}
While for $N=2M+1$ odd, we have 
$ y_{M+1}(\alpha) = \left ( 2M+\frac{3}{2} \right ) \pi$. In both cases, $y_k(\alpha) \approx (2k-1) \pi$ for $k > M+1$, as the saddle points above the $(M+1)$-th point are to occur around intersection of the local minima of $H^{\alpha}_1(x,y)$ with the axis, corresponding to local minima of $cos(y)$ (as seen in Fig. 2). Note also that by direct substitution, for $\alpha =  \frac{2}{(2N-1) \pi}$, we have 
\begin{equation}
\label{eq:zero} 
w_N \left (\alpha \right )=  (2N-1) \pi i.
\end{equation}

It should be remarked that approximations of $x_k(\alpha)$ can be obtained via the Lagrange inversion formula (see, for instance, 2.2 of \cite{DB}), by evaluating the solution of the equation 
$
H^{\alpha}_1(x,\widetilde{y}_k(\alpha))=0$, where $\widetilde{y}_k(\alpha)$ is the approximated value of $y_k(\alpha)$ discussed above. In particular, the growth of $x_k$ is $O(log(k))$ as $k \rightarrow \infty$. The following result could be viewed as a manifestation of the "exponential cancellation" phenomena for $\widetilde{A}(u,s)$.: 
\begin{cor} 
\label{cor:2.3}
For $s=\sigma+it$ and $u= \alpha \cdot t$ the following holds
\begin{equation}
\widetilde{A}(u,s)= O \left ( e^{ -\omega(\alpha) t} \right ), 
\end{equation}
where 
\begin{equation} 
\omega(\alpha):=- \alpha \cdot log \abs*{\frac{1-e^{-w_{k(\alpha)}(\alpha)}}{2}} + \theta_{k(\alpha)}(\alpha)-\frac{\pi}{2}, 
\end{equation}
and $w_k(\alpha):=r_k(\alpha)e^{i \theta_k(\alpha)}$ is the polar presentation of the $k(\alpha)$-th saddle-point.
\end{cor}

\hspace{-0.6cm} \bf Proof: \rm By (\ref{eq:A}) and Proposition \ref{prop:2.2}
\begin{multline} 
\label{eq:42}
 \vert 
\widetilde{A}(u,s)\vert  = O \left ( \sqrt{\frac{2 \pi}{t \vert f''(w_{k(\alpha)}(\alpha) ;\alpha ) \vert }} \cdot \frac{e^{-w_{k(\alpha)}(\alpha)}}{2 (1-2^{1-s}) \Gamma(s)}  \left (\frac{ 1-e^{-w_{k(\alpha)}(\alpha)}}{2} \right )^{\alpha t} w_{k(\alpha)}(\alpha)^{\sigma +it -1} \right )= \\ =O \left ( \sqrt{\frac{1}{t  }} \cdot \left (\frac{ 1-e^{-w_{k(\alpha)}(\alpha)}}{2} \right )^{\alpha t} \frac{w_{k(\alpha)}(\alpha)^{\sigma +it -1}}{\Gamma( \sigma+it)} \right )= O \left ( \left (\frac{ 1-e^{-w_{k(\alpha)}(\alpha)}}{2} \right )^{\alpha t} \left (\frac{w_{k(\alpha)}(\alpha)}{\sigma+it} \right)^{\sigma +it }\right ).
\end{multline} In the second step we use $\Gamma(\sigma+it)=O \left ( (\sigma+it)^{\sigma+it-\frac{1}{2}} \right )$, which follows from Stirling's formula. $\square$

\bigskip

In particular, in the setting of Theorem \ref{thm:main-theorem} of the introduction, where $\vert t \vert = 2 \pi x y$, we set $
\omega(x,y):=\omega(\alpha)$, with $\alpha = \frac{1}{2\pi y}$. 

Let us conclude this section by noting that, for a given $\alpha>0$, the  solutions $w_k(\alpha)$ can be evaluated numerically to any level of percision. For instance, by the Luck and Stevens formula, see \cite{LS,LZC}. For instance, for $\alpha=\frac{1}{\pi}$, we have for the leading saddle point 
\begin{equation}
\label{eq:Spoint}
w_2 \left ( \alpha \right ) \approx 0.68154 + 9.31481 i,
\end{equation} 
and hence $\omega(\alpha) \approx 0.017728$. The following Fig.3 illustrates Corollary \ref{cor:2.3} in this case: 
\begin{figure}[ht!]
	\centering
	\includegraphics[scale=0.35]{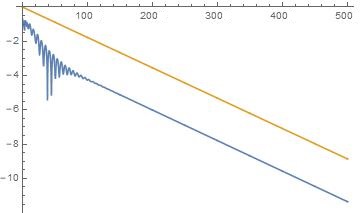}
	\caption{Graphs of $-\omega(\alpha)t$ (brown) and $log \left ( \widetilde{A}(\alpha \cdot t,0.05+i
	t) \right )$ (blue) for $\alpha=\frac{1}{\pi}, \omega(\alpha) = 0.017728$ and  $0 \leq y \leq 500$. }
\end{figure} 

 \section{The Asymptotics for the Integral of \texorpdfstring{$\widetilde{A}(u,s)$}{$I_1(x,s)$}}

In this section we consider the asymptotic behavior of the integral 
$
\widetilde{I}_1(x,s):=\int_x^{\infty} \widetilde{A}(u,s) du$.
We have:
\begin{theorem}
\label{thm:first-integral}
Let $s=\sigma+it,x = \alpha \cdot t$ and $N \in \mathbb{N}$. In $ \alpha \in \left [ \frac{2 }{ (N+1) \pi}, \frac{2 }{N \pi } \right ]$ the following holds
\begin{equation} 
\widetilde{I}_1(x,s):= \frac{\chi(s)}{2^{s-1}-1} \left (\sum_{k=1}^{N} (2k-1)^{s-1} \right ) +O \left (e^{-\omega(\alpha) t} \right ).
\end{equation}
\end{theorem}

\hspace{-0.6cm} \bf Proof: \rm
Integrating (\ref{eq:GL}) we obtain  
\begin{multline}
\label{eq:Int1}
\widetilde{I}_1(x,s)=  \frac{1}{2(1-2^{1-s})\Gamma(s)}  \int_x^{\infty}  \left ( \int_{0}^{\infty}  
e^{-w} \left ( \frac{ 1- e^{-w}}{2} \right)^u w^{s-1} dw \right ) du
= \\ =
 \frac{1}{2(1-2^{1-s})\Gamma(s)} \int_0^{\infty}    
e^{-w} \left ( \int_x^{\infty} \left ( \frac{ 1- e^{-w}}{2} \right)^u du \right ) w^{s-1} dw  = \\ =\frac{1}{2(1-2^{1-s})\Gamma(s)}  \int_{0}^{\infty}  \left [ \frac{  e^{-w}\left (\frac{1-e^{-w}}{2} \right )^x  }{ln \left ( \frac{1-e^{-w}}{2} \right ) } \right ] w^{s-1} dw. \end{multline}
As in (\ref{eq:29}), let us set \begin{equation}
\begin{array}{ccc} g(\sigma,w):=e^{-w} w^{\sigma -1} 
& ; & f(w; \alpha):=\alpha \cdot ln \left (1-e^{-w} \right ) + i \cdot ln(w), 
\end{array}
\end{equation}
and
\begin{equation}
h(w):= ln \left (\frac{1-e^{-w}}{2} \right ). 
\end{equation} 
For $\alpha \in \mathbb{R}$ let us set consider the $1$-parametric family of integrals:
\begin{equation}
\widetilde{I}_1(t,\sigma ; \alpha):= \frac{1}{2(1-2^{1-s})\Gamma(s)}  \int_{0}^{\infty} \frac{g(\sigma,w)}{h(w)} e^{t f(w; \alpha)}  dw. \end{equation}
The integrand of $\widetilde{I}_1(t,\sigma; \alpha)$ admits poles\footnote{ Historically, the asymptotic evaluation of integrals with poles was first studied in relation to the investigation of propagation of radio waves over a plane
earth, see \cite{So}. The question attracted the attention of mathematicians such as Van der Werden \cite{vanderw} and H. Weyl \cite{We}. We refer the reader also to the book by Felsen and Marcuvitz \cite{FM} which contains a detailed discussion on saddle-point techniques in the case of poles and branches.} at $(2k-1) \pi i$, for $k \in \mathbb{Z}$, coming from the zeros of $h(w)$. In particular, the saddle point technique requires us to deform the original contour of integration $[0,\infty)$, to a new contour, $C_{\alpha} \subset \mathbb{C}$, passing through the main contributing saddle point. Considering $C_{\alpha}$ as a continuous $1$-parametric family, the contour $C_{\alpha}$ is required to cross the first $N$ poles. Figure 3 shows a schematic illustration: 
\begin{figure}[ht!]
	\centering
	\includegraphics[scale=0.175]{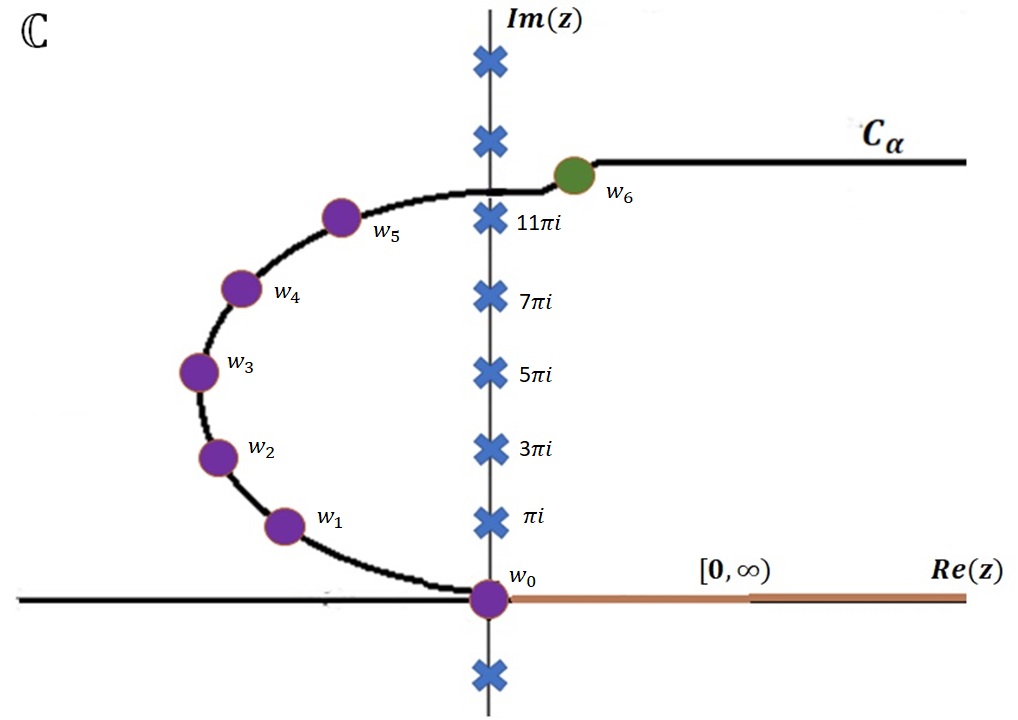}
	\caption{$C_{\alpha} \subset \mathbb{C}$ with singularities (blue) and saddle points (purple$\setminus$green). }
\end{figure}

By Cauchy's theorem we have 
\begin{equation}
\int_{0}^{\infty} \frac{g(\sigma,w)}{h(w)} e^{t f(w; \alpha)}  dw =\int_{C_{\alpha}} \frac{g(\sigma,w)}{h(w)} e^{t f(w; \alpha)}  dw +\sum_{k=1}^N Res \left (  \frac{g(\sigma,w)}{h(w)} e^{t f(w; \alpha)} ; (2k-1) \pi i \right ) 
\end{equation}
 In particular, each pole contributes the following residue 
\begin{multline} 
\label{eq:Residue2}
  Res \left (  \frac{g(\sigma,w)}{h(w)} e^{t f(w; \alpha)} ; (2k-1) \pi i \right ) = \\ =- 2 \pi i \cdot Res \left ( \frac{1}{h(w)}; (2k-1)\pi i \right )   g(\sigma, (2k-1)\pi i) e^{t f((2k-1) \pi i; \alpha)  }= \\=- 2 \pi i \cdot (-2) \cdot 
e^{-\pi i} \left ( \frac{ 1- e^{-\pi i}}{2} \right)^u ((2k-1) \pi i)^{s-1}  =-4 (2k-1)^{s-1} \pi^s e^{ \frac{-\pi i s}{2}}.
\end{multline}
The functional equation implies $\chi(s) = (\chi(1-s))^{-1}$. Hence, one can express 
$
\chi(z)=\frac{2^{s-1} \pi^{s}}{ cos \left ( \frac{s \pi}{2} \right ) \Gamma(s)  }$.
Clearly, 
\begin{equation} 
\frac{1}{ cos \left ( \frac{s \pi}{2} \right ) }=2 e^{-\frac{\pi i s}{2}}+O \left ( e^{-\frac{3 \pi s}{2} } \right ).
\end{equation}
Thus, 
\begin{multline} 
\label{eq:Residue}
\frac{1}{2(1-2^{1-s}) \Gamma(s)}  Res \left (  \frac{g(\sigma,w)}{h(w)} e^{t f(w; \alpha)} ; (2k-1) \pi i \right )  =-\frac{4 (2k-1)^{s-1} \pi^s e^{ \frac{-\pi i s}{2}}}{2(1-2^{1-s}) \Gamma(s)}= \\=\frac{\chi(s)}{1-2^{s-1}} \left ( (2k-1)^{s-1} +O \left ( e^{-\frac{3 \pi s}{2}} \right ) \right ).
\end{multline}
Finally, the main saddle point contributing to the integral over $C_{\alpha}$ is the $k(\alpha)$-th, as required. $\square$

\section{The Asymptotics for the Integral of $I_2(x,s)$}

In this section we consider the asymptotic behavior of the integral 
\begin{equation} 
\label{eq:der} \widetilde{I}_2(x,s):= \int_x^{\infty} \psi(u) \left ( \frac{\partial}{\partial u} \widetilde{A}(u,s) \right ) du,
\end{equation} 
where $
\psi(u):=u-[u]-\frac{1}{2}$.
We have:

\begin{theorem}
\label{thm:int2}
Let $s=\sigma+it$ and $x=\alpha \cdot t$. Then the following holds
\begin{equation} 
\widetilde{I}_2(x,s):= O \left (e^{-\omega(\alpha) t} \right ).
\end{equation}

\end{theorem}

\hspace{-0.6cm} \bf Proof: \rm As in the classical proof of Hardy and Littlewood in \cite{HL}, consider the Fourier series expansion of $\psi(u)$, which is given by 
\begin{equation}
 \psi(u)= - \sum_{n=1}^{\infty} \frac{1}{n \pi} sin(2n \pi u).
 \end{equation}
 Hence, we have 
 \begin{equation}
 \label{eq:60} \widetilde{I}_2(x,s):=- \sum_{n=1}^{\infty} \frac{1}{n \pi} \widetilde{I}_2^n(x,s),
\end{equation}
where
\begin{equation} \widetilde{I}_2^n(x,s):= \int_x^{\infty} sin(2n \pi u) \left ( \frac{\partial}{\partial u} \widetilde{A}(u,s) \right ) du.
\end{equation} 
 Let us express $
2i \widetilde{I}_2^n(x,s)= \widetilde{J}^+_n(x,s)+\widetilde{J}^-_n(x,s)$, where 
\begin{equation} 
 \widetilde{J}^{\pm}_n(x,s):=  \int_x^{\infty} e^{\pm 2n \pi i u} \left ( \frac{\partial}{\partial u} \widetilde{A}(u,s) \right ) du.
\end{equation}
By direct expansion, we have  
\begin{equation}
\widetilde{J}_n^{\pm}(x,s)= \frac{1}{2(1-2^{1-s})}  \int_0^{\infty} \left ( \frac{ ln \left ( \frac{1-e^{-w}}{2} \right )  .\left ( \frac{ 1- e^{-w}}{2} \right)^x  e^{-w \pm 2 n \pi i x}}{ ln \left ( \frac{1-e^{-w}}{2} \right ) \pm 2n\pi i }   \right )  w^{s-1} dw 
\end{equation}
For $s=\sigma+it$ and $u = \alpha \cdot t$ with $t \rightarrow \infty$ set 
\begin{equation}
\widetilde{J}_n^{\pm}(\sigma,t ; \alpha)=\frac{\widetilde{k}^{\pm}_n(\sigma,t ; \alpha)}{2^{x+1}(1-2^{1-s}) \Gamma(s) },
\end{equation}
where 
\begin{equation}
\label{eq:II}
\widetilde{k}^{\pm}_n(\sigma,t ;\alpha ) := \int_0^{\infty} \widetilde{g}_n^{\pm}(\sigma,w) e^{t \widetilde{f}_n^{\pm}(w ; \alpha )} dw,
\end{equation} with 
\begin{equation}
\begin{array}{ccc} \widetilde{g}^{\pm}_n(\sigma,w):=\frac{ln \left ( \frac{1-e^{-w}}{2} \right ) e^{-w} w^{\sigma -1}}{ln \left ( \frac{1-e^{-w}}{2} \right ) \pm 2n \pi i} 
& ; & \widetilde{f}^{\pm}_n(w ; \alpha):= \alpha \cdot ln \left (1-e^{-w} \right ) + i \cdot \left ( ln(w) \pm 2n \right ). 
\end{array}
\end{equation}
On the one hand, the integrand, in this case, has no poles since, as for $n \neq 0$, there are no solutions of 
$ ln \left ( \frac{1-e^{-w}}{2} \right ) = \pm 2n \pi i$,
for the principle branch of the complex logarithm (in contrast to the case of $\widetilde{I}_1(x,s)$ of the previous Section 3). 
On the other hand, note that for $k \in \mathbb{Z}$, the integrand does admit branch cuts, starting at $2k \pi i$ and extending to infinity from the left. However, these branch cuts do not contribute to the asymptotic behavior of $\widetilde{J}^{\pm}_n(\sigma,t ;\alpha)$, since it is not necessary to cross the cuts when deforming $[0,\infty)$ to the saddle-point contour. Figure 4 shows a schematic illustration of the saddle point contour, avoiding crossing of the branch cuts. 
\begin{figure}[ht!]
	\centering
	\includegraphics[scale=0.18]{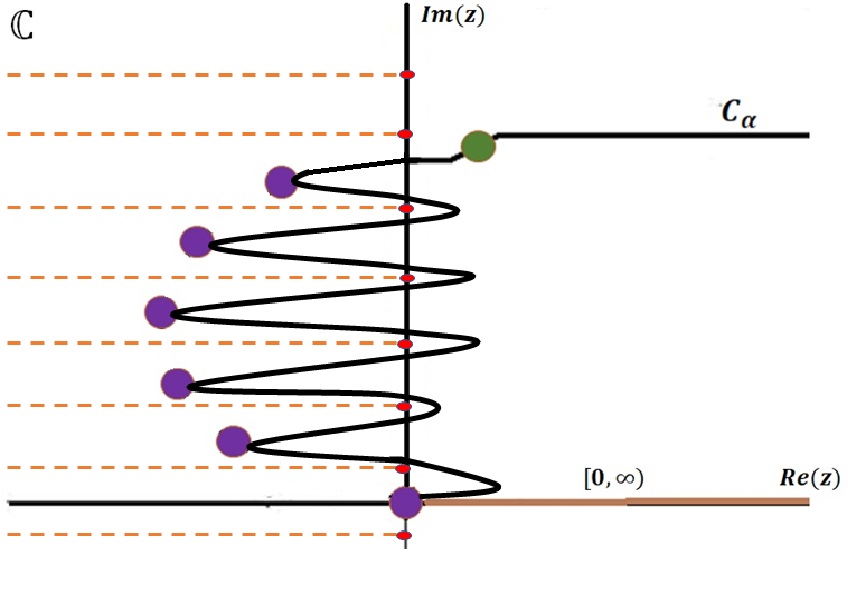}
	\caption{Integration contour with branch cuts (red) and saddle points (purple$\setminus$green). }
\end{figure}

Hence, only saddle points contribute to the asymptotic properties of the integral $\widetilde{k}^{\pm}_n(\sigma,t,\alpha)$. The saddle points are given by the solutions of 
\begin{equation}
 \frac{\partial \widetilde{f}^{\pm}_n}{\partial w} (w ; \alpha) = \frac{ \alpha \cdot e^{-w}}{1-e^{-w}}+\frac{i}{w}=0. 
\end{equation} which is similar to the saddle-point equation (\ref{eq:equ30}) of Section 2. We, hence, have 
\begin{multline} 
\widetilde{J}^{\pm}_n(x,s)= O \left (  \sqrt{\frac{2 \pi}{t \vert f''(w_k(\alpha); \alpha) \vert}}    \frac{ ln \left ( \frac{1-e^{-w_{k(\alpha)}(\alpha)}}{2} \right )  \left ( \frac{ 1- e^{-w_{k(\alpha)}(\alpha)}}{2} \right)^x  e^{-w_{k(\alpha)}(\alpha) \pm 2 n \pi i x}}{ 2(1-2^{1-s}) \left (ln \left ( \frac{1-e^{-w_{k(\alpha)}(\alpha)}}{2} \right ) \pm 2n\pi i \right )\Gamma(s) }   w_{k(\alpha)}(\alpha)^{s-1}  \right ) =\\ = O \left ( \frac{1}{ n \cdot \Gamma(s) }  \sqrt{\frac{1}{t }}    \left ( \frac{ 1- e^{-w_{k(\alpha)}(\alpha)}}{2} \right)^x   e^{-w_{k(\alpha)}(\alpha)} w_{k(\alpha)}(\alpha)^{s-1}  \right )=O \left ( \frac{  e^{-\omega(\alpha) t}}{n} \right ).
\end{multline} 
By definition, also $\widetilde{I}_2^n(x,s)=O \left ( \frac{  e^{-\omega(\alpha) t}}{n} \right )$. Hence, by (\ref{eq:60}) we have 
\begin{equation} 
\widetilde{I}_2(x,s)=O \left (\sum_{n=1}^{\infty} \frac{  e^{-\omega(\alpha) t}}{n^2} \right )=O \left (e^{-\omega(\alpha) t} \right ),
\end{equation}
as required. $\square$

 \section{Summary and Concluding Remarks}
 
In this work we have proved a new approximate functional equation, whose error term is of exponentional rate of decay (\ref{eq:NEW-AFE}). Our proof, as the proof of Hardy and Littelwood in the classical case \cite{HL}, relied on the Euler-Maclaurain formula (\ref{eq:EM}). However, it is interesting to draw a few distinctions and similarities between the two settings.

First, let us note that, in the classical case, a direct application of the Euler-Maclaurin formula, as in the proof of Theorem \ref{thm:main-theorem}, does not lead to a proof of the full approximate functional equation (\ref{eq:HL-AFE}). Instead, this approach, in the classical case, leads only to the so-called "simplest approximate functional equation"\footnote{ See Theorem 4.11 of \cite{T} and Theorem 1.8 of \cite{I}}  which states 
\begin{equation}
\label{AFE-simple}
\zeta(s)= \sum_{n \leq x} n^{-s} +\frac{x^{1-s}}{s-1}+O(x^{-\sigma}),
\end{equation} 
for $\vert t \vert \leq \pi x$, see Lemma 2 of \cite{HL}. In fact, the main bulk of \cite{HL}, as well as \cite{HL3}, involves the development of ingenious analytic methods to extend the proof of the "simple" equation to a proof of the full classical approximate functional equation (\ref{eq:HL-AFE}).  

Concretely, as mentioned in the introduction, in the classical case, one applies the Euler-Maclaurain formula (\ref{eq:EM}) to the function $A(u,s)=u^{-s}$. As a result, the corresponding first integral, in this case, is given by $
I_1(x,s)= \int_x^{\infty}A(u,s) du = -\frac{x^{1-s}}{1-s}$ in the domain $\sigma>1$. In particular, the summation formula can be written in the following form 
\begin{equation}
\label{eq:simple2}
\zeta(s)-\sum_{n \leq x} A(u,s) =-\frac{x^{1-s}}{1-s}+ I_2(x,s) +O(x^{-s}),
\end{equation} where $
I_2(x,s):=\int_x^{\infty} \psi(u) \left ( \frac{\partial}{\partial u} A(u,s) \right ) du$. Note that (\ref{eq:simple2}) still holds, by analytic continuation, for $\sigma >-2$. The proof of the "simple approximate equation" is thus based on showing that 
$
I_2(x,s)=O(x^{-s})$, when $\vert t \vert \leq \pi x$. However, complications arise when considering the values in the range $\pi x \leq \vert t \vert$. In this case, the difference between the value of $\zeta(s)$ and the $\left [x \right]$-th partial sum of the series $
D(x,s):=\zeta(s)-\sum_{n \leq x } A(u,s)$ starts to largely deviate  from the first integral $I_1(x,s)=-\frac{x^{1-s}}{1-s}$. Hence, the power of the full (classical) approximate functional equation is required. 

In contrast, in the Hasse-Sondow setting our Theorem \ref{thm:int2} shows that the analog second integral is of exponential rate of decay 
$
\widetilde{I}_2(x,s)=O(e^{-\omega(x,y)})$, with $\vert t \vert = \pi x y$ for any $y \geq 1$. Correspondingly, Theorem \ref{thm:first-integral} shows that the first integral 
\begin{equation}
\widetilde{I}_1(x,s):= \frac{\chi(s)}{2^{s-1}-1} \left (\sum_{k=1}^{N} (2k-1)^{s-1} \right ) +O \left (e^{-\omega(x,y) t} \right )
\end{equation} serves as an approximation of 
$\widetilde{D}(x,s):=\zeta(s)-\sum_{n \leq x} \widetilde{A}(u,s)$,
 for any value of $x$. In particular, this sensitivity to the values of $x$, comes from the contribution of the added poles of the integrand, as seen in the proof of Theorem \ref{thm:first-integral}. These unique features of the Hasse-Sondow setting, enable us to actually obtain the proof of Theorem \ref{thm:main-theorem}, which is a full approximate functional equation, along the lines of the proof of the "simple approximate functional equation", in the classical case. However, with the trade-off that more involved saddle-point techniques are required to be applied.  

Finally, it is interesting to note a possible relation of the above results to the field of quantum chaos. In \cite{BK} Berry and Keating conjecture that the imaginary values $t_n$ of zeros of zeta could, in theory, be realized as the eigenvalues of an Hermitean operator $H$, obtained by quantizing a certain classical dynamical system with Hamiltonian $H_{cl}$. In this framework, the special form of the Riemann-Siegel formula plays a substantial role. Concretely, it is suggested that, under the assumption that $H_{cl}$ exists, the main sum of the Riemann-Siegel formula arises from contributions to the corresponding quantum spectral determinant, coming from periodic orbits of the Hamiltonian with real energy levels, while the correction terms of the Riemann-Siegel formula are expected to arise from periodic orbits of complex energy levels. Based on such insights, as well as extensive numerical evidence, Berry \cite{Be} and Berry-Keating \cite{BK} conjecture that the optimal accuracy of the Riemann-Siegel formula is $O(e^{- \pi t})$, and that this  exponential error should arise from the existence of such complex periodic orbits, see also \cite{Be2}. In the classical Riemann-Siegel setting such an expected optimal error $O(e^{-\pi t})$
is rather elusive, as when expanded to a given fixed order the error of the Riemann-Siegel formula is of bigger magnitude than exponential. However, in our setting of Theorem \ref{thm:main-theorem} the resulting exponential term arises directly from the leading solution of the transcendental equation $e^w-\alpha \cdot i w -1=0$, as well as weaker exponential terms arising from the other solutions. It is thus interesting to suggest that saddle points, which are the solutions of this transcendental equation, actually correspond to the complex periodic orbits expected by the Berry-Keating framework, in a similar manner to the way that pole contributions correspond to real periodic orbits.


\begin{thebibliography}{10}

\bibitem{Be} M. V. ~Berry. \newblock The Riemann-Siegel Expansion for the Zeta Function: High Orders and Remainders. \newblock Proceedings: Mathematical and Physical Sciences, vol. 450, no. 1939, 1995, pp. 439--462. 

\bibitem{Be2} M. V. ~Berry. \newblock Riemann's Saddle-point Method and the Riemann-Siegel Formula. \newblock The Legacy of Bernhard Riemann
After One Hundred and Fifty Years
ALM35, pp. 69–78, Higher Education Press
and International Press
Beijing–Boston.

\bibitem{BK} M. V. Berry and J. P. Keating. \newblock The Riemann Zeros and Eigenvalue Asymptotics. \newblock SIAM Review 41(2), 236-266, November 2001.

\bibitem{BH} N.~Bleistein and R. A. Handelsman. \newblock Asymptotic Expansions of Integrals. \newblock Holt, Rinehart \& Winston, New York, 1977. 

\bibitem{DB} N. G. de Bruijn. \newblock Asymptotic methods in analysis. \newblock Bibliotheca Mathematica. Vol. 4, North-Holland Publishing Co., Amsterdam; P. Noordhoff Ltd., Groningen; Interscience Publishers Inc., New York, 1958. 


\bibitem{E} H. M. Edwards. \newblock Riemann's zeta function. \newblock Pure and Applied Mathematics, 1974, 58, New York-London: Academic Press.

\bibitem{FM}  L. B. Felsen, N. Marcuvitz. 
\newblock Radiation and Scattering of Waves. \newblock Prentice-Hall, Microwaves and Fields Series, 1973.  

\bibitem{FS} P.~Flajolet and R.~Sedgewick. 
\newblock Mellin transforms and asymptotics: Finite differences and Rice's integrals. \newblock Theoretical Computer Science 144 (1995) pp 101--124.	

\bibitem{GN}
M. L. Glasser and I. Nagy. \newblock
Moments of powers of the Hulth{\'e}n density. \newblock Journal of Mathematical Chemistry, 2012,
50, 1707-1710.

\bibitem{GL} S.~Goyal and R. K. Laddha. 
\newblock On the generalized Riemann zeta funcion and the generalized Lambert transform. \newblock Ganita Sandesh, 11, 99-108, 1997.



\bibitem{HL} G. H.~Hardy and J. E.~Littlewood. \newblock The zeros of Riemann's zeta function on the critical line. \newblock Math. Z. 10, 283--317 (1921).
	
\bibitem{HL2} G. H.~Hardy and J. E.~Littlewood. \newblock The approximate functional equation in the theory of the zeta function, with an application to the divisor-problems of Dirichlet and Piltz. \newblock Proceedings of the London Mathematical Society, Volume s2-21, Issue 1, 1923, Pages 39--74.


\bibitem{HL3} G. H.~Hardy and J. E.~Littlewood. \newblock The approximate functional equations for $\zeta(s)$ and $\zeta^2(s)$. \newblock Proceedings of the London Mathematical Society, S2-29 (1), 81-97, 1929. 

\bibitem{H} H.~Hasse. \newblock Ein Summierungsverfahren fur die Riemannsche $\zeta$-Reihe. \newblock Math. Z. 32: 458--464 (1930).

\bibitem{I} A. Ivic. \newblock The Riemann Zeta Function. John Wiley \& Sons, 1985. 

\bibitem{K} D. E. Knuth. \newblock The art of computer programming, Vol. 3: Sorting and Searching. \newblock Addison-Wesley Reading, MA, 1973. 

\bibitem{LS} R. Luck and J.W. Stevens. \newblock Explicit solutions for transcendental
equations. \newblock SIAM Review, vol. 44, no. 2, pp. 227--233, 2002.

\bibitem{LZC} R. Luck, G. J. Zdaniuk, and H. Cho. \newblock An Efficient Method to Find Solutions for Transcendental Equations with Several Roots. \newblock International Journal of Engineering Mathematics, vol. 2015, Article ID 523043, 4 pages, 2015.

\bibitem{NC} S. Nadarajah and J. Chu. \newblock On moments of powers of the Hulth{\'e}n density. \newblock Journal of Mathematical Chemistry, 2016, 55, 911--913.

\bibitem{N} N. E. N{\o}rlund. \newblock Vorlesungen uber Differenzenrechnung. \newblock (1954) Chelsea Publishing Company, New York.


	\bibitem{SE} J. Sondow. \newblock Analytic continuation of Riemann's zeta function and values at negative integers via Euler's transformation of 
	series. \newblock Proceedings of the American mathematical society, Volume 120, Issue 2, 1994, 421--425.
	
	\bibitem{Si} C. L. Siegel. \newblock Über Riemanns Nachlaß zur analytischen Zahlentheorie. \newblock Quellen Studien zur Geschichte der Math. Astron. und Phys. Abt. B: Studien 2: 45--80, 1932. (Also in Gesammelte Abhandlungen, Vol. 1. Berlin: Springer-Verlag, 1966).
	
	\bibitem{So} A. Sommerfeld. \newblock {\"U}ber die Ausbreitung der Wellen in der drahtlosen Telegraphie. \newblock Ann. Physik .8 (1909) 665--736. 
	
	\bibitem{T} E. C. Titchmarsh and D. R. Heath-Brown (ed.). \newblock The Theory of the Riemann Zeta Function, 1984, (2nd rev. ed.). Oxford University Press.
	
	\bibitem{vanderw} B. L. van der Waerden. \newblock On the method of saddle points. \newblock Applied Scientific Research, Section B, December 1952, Volume 2, Issue 1, 33--45.
	
	\bibitem{We} H.~Weyl. \newblock Ausbreitung elektromagnetischer Wellen {\"u}ber einem ebenen Leiter. \newblock Ann. Physik 60 (1919) 481-500.
	
	\bibitem{Wo} R.~Wong. \newblock  Asymptotic approximations of integrals. \newblock  Classics in Applied Mathematics, SIAM, 2001. 
	 
\end{thebibliography}
\end{document}